\documentclass[12pt,a4paper]{amsart}
\usepackage{amsmath}
\usepackage{amscd}
\usepackage{amssymb}
\numberwithin{equation}{section}

\input{prepictex}
\input{pictex}
\input{postpictex}
\chardef\bslash=`\\ 

\makeatletter
\def\verbatim{\interlinepenalty\@M \@verbatim
  \leftskip\@totalleftmargin\advance\leftskip2pc
  \frenchspacing\@vobeyspaces \@xverbatim}
\makeatother
\newtheorem{theorem}{Theorem}[section]

\newtheorem{lemma}[theorem]{Lemma}
\newtheorem{proposition}[theorem]{Proposition}

\theoremstyle{definition}
\newtheorem{definition}[theorem]{Definition}
\newtheorem{remark}[theorem]{Remark}

\newtheorem{example}[theorem]{Example}

\theoremstyle{main}
\newtheorem{main}{Theorem}
\theoremstyle{maincorollary}
\newtheorem{maincorollary}{Corollary}

\newcounter{picture}

\DeclareMathOperator{\covol}{covol}

\newcommand{\FF}{{\mathbb F}}

\newcommand{\PP}{{\mathbb P}}
\newcommand{\QQ}{{\mathbb Q}}
\newcommand{\RR}{{\mathbb R}}

\newcommand{\ZZ}{{\mathbb Z}}
\newcommand{\cA}{{\mathcal A}}

\newcommand{\cO}{{\mathcal O}}

\newcommand{\D}{{\Delta}}
\newcommand{\bD}{{\partial \D}}

\newcommand{\e}{{\varepsilon}}
\newcommand{\G}{{\Gamma}}
\newcommand{\Om}{{\Omega}}
\newcommand{\g}{{\gamma}}
\newcommand{\s}{{\sigma}}
\newcommand{\w}{{\omega}}

\newcommand{\fL}{{\mathfrak L}}

\newcommand{\od}{{\overline d}}

\newcommand{\ove}{{\overline e}}

\newcommand{\zomatrix}{{$\{0,1\}$-matrix}}

\newcommand{\id}{{\bf 1}}
\newcommand{\tl}{{{\text{\rm{III}}}_{\lambda}}}

\newcommand{\tone}{{{\text{\rm{III}}}_{1}}}
\newcommand{\tq}{{{\text{\rm{III}}}_{1/q}}}
\newcommand{\tqq}{{{\text{\rm{III}}}_{1/{q^2}}}}
\newcommand{\aut}{{\text{\rm Aut}}}

\newcommand{\PGL}{{\text{\rm{PGL}}}}
\newcommand{\PSL}{{\text{\rm{PSL}}}}

\newcommand{\Li}{{\text{\rm{L}}}^{\infty}}

\begin{document}

\title[]{Boundary operator algebras for free uniform tree lattices}

\author{Guyan Robertson}
\address{School of Mathematics and Statistics, University of Newcastle, NE1 7RU, U.K.}
\email{a.g.robertson@newcastle.ac.uk}
\date{\today}          
\thanks{This research was supported by the Australian Research Council.} 

\begin{abstract}
Let $X$ be a finite connected graph, each of whose vertices has degree at least three. The fundamental group $\G$ of $X$ is a free group and  acts on the universal covering tree $\D$ and on its boundary $\bD$, endowed with a natural topology and Borel measure. The crossed product $C^*$-algebra $C(\bD) \rtimes \G$ depends only on the rank of $\G$ and is a Cuntz-Krieger algebra whose structure is explicitly determined. The crossed product von Neumann algebra does not possess this rigidity. If $X$ is homogeneous of degree $q+1$ then the von Neumann algebra $\Li(\bD)\rtimes \G$ is the hyperfinite factor of type  $\tl$ where $\lambda=1/{q^2}$ if $X$ is bipartite, and $\lambda=1/{q}$ otherwise. 
\end{abstract}      

\maketitle

\section*{Introduction}

Let $\D$ be a locally finite tree whose automorphism group $\aut(\D)$ is equipped with the compact open topology. Let $\G$ be a discrete subgroup of $\aut(\D)$ which acts freely on $\D$. That is, no element $g\in \G-\{1\}$ stabilizes any vertex or geometric edge of $\D$. Assume furthermore that $\G$ acts cocompactly on $\D$, so that the quotient $\G\backslash\D$ is a finite graph. Then $\G$ is a finitely generated free group and is referred to as a free uniform tree lattice. 

Conversely, if $X$ is a finite connected graph and $\G$ is the fundamental group of $X$, then $\G$ is a finitely generated free group and acts freely and cocompactly on the universal covering tree $\D$.

It is fruitful to think of the tree $\D$ as a combinatorial analogue of the Poincar\'e disc and $\G$ as an analogue of a Fuchsian group.
The group $\G$ is the free group on $\g $ generators, where $\g =1-\chi(\G\backslash\D)$ and $\chi(\G\backslash\D)$ is the Euler-Poincar\'e characteristic of the quotient graph.
Let $S$ be a free set of generators for $\G$. 

Define a \zomatrix \  $A$ of order $2\g$, with entries indexed by elements of $S\cup S^{-1}$,  by
\begin{equation}
A(x,y) = 
\begin{cases}
    1&  \text{if $y\ne x^{-1}$},\\
    0&  \text{if $y= x^{-1}$.}
\end{cases}
\end{equation}
Notice that the matrix $A$ depends only on the rank of the free group $\G$.

The boundary $\bD$ of the tree $\D$ is the set of equivalence classes of infinite semi-geodesics in $\D$, where equivalent semi-geodesics contain a common sub-semi-geodesic. There is a natural compact totally disconnected topology on $\bD$ \cite[I.2.2]{ser}.
Denote by $C(\bD)$ the algebra of continuous complex valued functions on $\bD$. The full crossed product algebra $C(\bD) \rtimes \G$ is the universal $C^*$-algebra  generated by the commutative $C^*$-algebra $C(\bD)$ and the image of a unitary
representation $\pi$ of $\Gamma$, satisfying the covariance relation $f(g^{-1}\w) = \pi(g)\cdot f \cdot \pi(g)^{-1}(\w)$ for $f \in C(\bD)$, $g \in \Gamma$ and $\w\in\bD$ \cite[Chapter 7]{ped}.

\begin{main}\label{main}
Let $\D$ be a locally finite tree whose vertices all have degree at least three.
Let $\G$ be a free uniform lattice in $\aut(\D)$.  Then the boundary $C^*$-algebra 
$\cA(\G)=C(\bD) \rtimes \G$ depends only on the rank of $\G$, and $\G$ is itself determined by $K_0(\cA(\G))$. More precisely, 
\begin{itemize}
\item[(1)]
$\cA(\G)$ is isomorphic to the simple Cuntz-Krieger algebra 
${\mathcal O}_A$ associated with the matrix $A$;
\item[(2)]
$K_0(\cA(\G))= \ZZ^\g \oplus \ZZ/(\g -1)\ZZ$ and the class of the identity $[\id]$ is the generator of the summand $\ZZ/(\g -1)\ZZ$. Moreover $K_1(\cA(\G))= \ZZ^\g$.
\end{itemize}
\end{main}

The algebra  $\cA(\G)$ satisfies the hypotheses of the classification theorem of \cite{k},\cite{ph}. Therefore the isomorphism class of the algebra $\cA(\G)$ is determined by its K-theory together with the class of the identity in $K_0$.
The fact that the class $[\id]$ in $K_0$ has order equal to $-\chi(\G\backslash\D)$ strengthens the result of \cite[Section 1]{rob} and provides an exact analogy with the Fuchsian case \cite{ad}.

Theorem \ref{main} will be proved in Lemmas \ref{mainlemma} and \ref{Ktheory} below.
The key point in the proof is that the Cuntz-Krieger algebra $\cO_A$ is defined uniquely, up to isomorphism, by a finite number of generators and relations \cite{ck}, and it is possible to identify these explicitly in $\cA(\G)$. The original motivation for this result was the paper of J. Spielberg \cite{sp}, which showed that if $\G$ acts freely and transitively on the tree $\D$ then $\cA(\G)$ is a Cuntz-Krieger algebra. Higher rank analogues were studied in \cite{rs2}.

 There is a natural Borel measure on  $\bD$ and one may also consider the crossed product von Neumann algebra $\Li(\bD)\rtimes \G$. This is the von Neumann algebra arising from the classical group measure space construction of Murray and von Neumann \cite{su}. In contrast to Theorem \ref{main}, the structure of this algebra depends on the tree $\D$ and on the action of $\G$.
For simplicity, only the case where $\D$ is a homogeneous tree is considered.

\begin{main}\label{main2}
Let $\D$ be a homogeneous tree of degree $q+1$, where $q\ge 1$, and let $\G$ be a free uniform lattice in $\aut(\D)$.  Then $\Li(\bD)\rtimes \G$ is the hyperfinite factor of type  $\tl$ where
\begin{equation*}
\lambda = 
\begin{cases}
1/{q^2} &  \text{if the graph $\G\backslash\D$ is bipartite}, \\    
1/q &  \text{otherwise}.
\end{cases}
\end{equation*}
\end {main}

Theorem \ref{main2} will be proved in Section \ref{measure section}.
The result could equally well have been stated as a classification of the measure theoretic boundary actions up to orbit equivalence \cite{ho}.
The analogous result for a Fuchsian group $\G$ acting on the circle is  that
$\Li(S^1)\rtimes \G$ is the hyperfinite factor of type $\tone$ \cite{spa}.  

The special case of Theorem \ref{main2} where $\G$ acts freely and transitively on the vertices of $\D$ was dealt with in \cite{rr}. In that case $q$ is odd, $\G$ is the free group of rank $\frac{q+1}{2}$, and  $\Li(\bD)\rtimes \G$ is the hyperfinite factor of type  $\tq$. 
We remark that R. Okayasu \cite{ok} constructs similar algebras in a different way, but does not explicitly compute the value of $\lambda$.

There is a type map~$\tau$ defined on the vertices of~$\D$ and taking values in $\ZZ/2\ZZ$, defined as follows.
Fix a vertex $v_0\in \D$ and let $\tau(v)= d(v_0,v) \pmod{2}$,
where $d(u,v)$ denotes the usual graph distance between vertices of the tree.
The type map is independent of $v_0$, up to addition of $1 \pmod{2}$. 
It therefore induces a canonical partition of the vertex set of $\D$ into two classes, so that two vertices are in the same class if and only if the distance between them is even. An automorphism $g\in\aut(\D)$ is said to be {\em type preserving} if, for every vertex $v$, $\tau(gv)= \tau(v)$.  The graph $\G\backslash\D$ is bipartite if and only if the action of $\G$ is type preserving.

Let $\FF$ be a nonarchimedean local field with residue field of order $q$.
The Bruhat-Tits building associated with $\PGL(2,\FF)$ is a regular tree $\D$ of degree $q+1$ whose boundary may be identified with the projective line $\PP_1(\FF)$.
If $\G$ is a torsion free lattice in $\PGL(2,\FF)$ then $\G$ is necessarily a free group of rank $\g \ge 2$, which acts freely and cocompactly on $\D$ \cite [Chapitres I.3.3, II.1.5]{ser}, and the results apply to the action of $\G$ on $\PP_1(\FF)$.

Let $\cO$ denote the valuation ring of $\FF$. Then $K=\PGL(2,\cO)$ is an open maximal compact subgroup of $\PGL(2,\FF)$ and the vertex set of $\D$ may be identified with the homogeneous space $\PGL(2,\FF)/K$. If the Haar measure $\mu$ on $\PGL(2,\FF)$ is normalized so that $\mu(K)=1$,
then the covolume $\covol(\G)$ is equal to the number of vertices of $X=\G\backslash\D$ and  $\g -1=\frac{(q-1)}{2}\covol(\G)$, (c.f. \cite[Chapitre II.1.5]{ser}).

The action of $\G$ on $\D$ is type preserving if and only if $\G$ is a subgroup of $\PSL(2,\FF)$.  Combining Theorem \ref{main} and Theorem \ref{main2}, in this special case, yields

\newpage

\begin{maincorollary}\label{maincorollary}
Let $\G$ be a torsion free lattice in $\PGL(2,\FF)$. Using the above notation, the boundary algebras are determined as follows. 
\begin{itemize}
\item[(1)] The $C^*$-algebra $\cA(\G)=C(\PP_1(\FF)) \rtimes \G$ is the unique Cuntz-Krieger algebra satisfying
\[
(K_0(\cA(\G)), [\id])= (\ZZ^\g \oplus \ZZ/(\g -1)\ZZ, (0,0,\dots, 0, 1)).
\]
\item[(2)] The von Neumann algebra $\Li(\PP_1(\FF))\rtimes \G$ is the hyperfinite factor of type 
$\tl$ where
\begin{equation*}
\lambda = 
\begin{cases}
1/{q^2} &  \text{if $\G\subset \PSL(2,\FF)$}, \\    
1/q &  \text{otherwise}.
\end{cases}
\end{equation*}
\end{itemize}
\end{maincorollary}

\bigskip

\section{The Cuntz-Krieger algebra}\label{section1}

Let $\D$ be a locally finite tree whose vertices all have degree at least three.
The results and terminology of \cite{ser} will be used extensively. 
The edges of $\D$ are directed edges and each geometric edge of $\D$ corresponds to two directed edges $d$ and $\od$.
Let $\D^0$ denote the set of vertices and $\D^1$  the set of directed edges of $\D$.

Suppose that $\G$ is a torsion free discrete group acting freely on $\D$ : that is no element $g\in \G-\{1\}$ stabilizes any vertex or geometric edge of $\D$. 
Then $\G$ is a free group \cite[I.3.3]{ser} and there is an orientation on the edges which is invariant under $\G$ \cite[I.3.1]{ser}. Choose such an orientation.  
This orientation consists of a partition $\D^1=\D^1_+\sqcup \overline{\D^1_+}$ and a bijective involution $d\mapsto \overline d : \D^1 \to \D^1$ which interchanges the two components of $\D^1$. Each directed edge $d$ has an
origin $o(d)\in \D^0$ and a terminal vertex $t(d)\in \D^0$ such that $o(\overline d)=t(d)$. 

Assume that $\G$ acts cocompactly on $\D$. This means that the quotient $\G\backslash\D$ is a finite connected graph with vertex set $V=\G\backslash\D^0$
and directed edge set $E=E_+\sqcup\overline{E_+}=\G \backslash\D^1_+\sqcup \G\backslash\overline{\D^1_+}$.  The Euler-Poincar\'e characteristic of the graph is $\chi(\G\backslash\D)=n_0-n_1$ where $n_0=\#(V)$ and $n_1=\#(E_+)$,
and $\G$ is the free group on $\g $ generators, where $\g =1-\chi(\G\backslash\D)$.

Choose a tree $T$ of representatives of $\D \pmod \G$; that is a lifting of a maximal tree of $\G\backslash\D$. The tree $T$ is finite, since $\G$ acts cocompactly on $\D$.
Let $S$ be the set of elements $x\in\G-\{1\}$ such that there exists an edge $e\in \D^1_+$ with $o(e)\in T$ and $t(e)\in xT$. Then $S$ is a free set of generators for the free group $\G$ \cite[I.3.3, Th\'eor\`eme $4^\prime$]{ser} and $\g =\# S$. It is clear that
$S^{-1}$ is the set of elements $x\in\G-\{1\}$ such that there exists an edge $e\in \D^1_-$ with $o(e)\in T$ and $t(e)\in xT$.
The map $g\mapsto gT$ is a bijection from $\G$ onto the set of $\G$ translates of the tree $T$ in $\D$, and these translates are pairwise disjoint \cite[I.3.3, Proof of Th\'eor\`eme $4^\prime$]{ser}. Moreover each vertex of $\D$ lies in precisely one of the sets $gT$.

The boundary $\bD$ of the tree $\D$ is the set of equivalence classes of infinite semi-geodesics in $\D$, where equivalent semi-geodesics agree except on finitely many edges. Also $\bD$ has a natural compact totally disconnected topology \cite[I.2.2]{ser}. The group $\G$ acts on $\bD$ and one can form the crossed product algebra $C(\bD) \rtimes \G$.
This is the universal $C^*$-algebra  generated by the commutative $C^*$-algebra $C(\bD)$ and the image of a unitary
representation $\pi$ of $\Gamma$, satisfying the covariance relation 
\begin{equation}\label{cov0}
f(g^{-1}\w) = \pi(g)\cdot f \cdot \pi(g)^{-1}(\w)
\end{equation} for
$f \in C(\bD)$, $g \in \Gamma$ and $\w\in\bD$ \cite{ped}.
This covariance relation implies that for each clopen set $E\subset\bD$ we have
\begin{equation}\label{cov}
\chi_{gE} = \pi(g)\cdot \chi_E\cdot  \pi(g)^{-1}.
\end{equation}
In this equation, $\chi_E$ is a continuous function and is regarded as an element of the crossed product algebra via the embedding $C(\bD)\subset C(\bD) \rtimes \G$. In the present setup the algebra $C(\bD) \rtimes \G$ is seen {\em a posteriori} to be simple. Therefore $C(\bD) \rtimes \G$ coincides with the reduced crossed product algebra \cite[7.7.4]{ped} and there is no need to distinguish between them notationally.

Fix a vertex $O\in \D$ with $O\in T$. Each $\w\in\bD$ has a unique representative semi-geodesic $[O,\w)$ with initial vertex $O$. A basic open neighbourhood of $\w\in \bD$ consists of those $\w^{\prime}\in \bD$ such that $[O,\w)\cap [O,\w^{\prime})\supset [O,v]$ for some fixed $v\in [O,\w)$.
If $g\in\G-\{1\}$, let $\Pi_g$ denote the set of all $\w \in\bD$ such that $[O,\w)$ meets the tree $gT$.
Note that $\Pi_g$ is clopen, since $T$ is finite. The characteristic function $p_g$ of the set $\Pi_g$ is continuous and so lies in $C(\bD)\subset C(\bD) \rtimes \G$. The identity element $\id$ of $C(\bD) \rtimes \G$ is the constant function defined by $\id(\w)=1,\, \w\in\bD$.

\begin{lemma}\label{L1}
If $x,y\in S\cup S^{-1}$ with $x\ne y^{-1}$ then 
\begin{itemize}
\item[(a)]
$\pi(x) p_{x^{-1}} \pi(x^{-1})= \id - p_x$ ;
\item[(b)]
$\pi(x) p_{y} \pi(x^{-1})= p_{xy}$.
\end{itemize}
\end{lemma}

\begin{figure}[htbp]
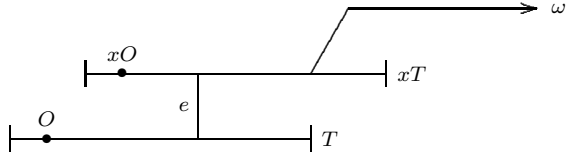
\label{fig1}
\hfil
\centerline{
\beginpicture
\setcoordinatesystem units <0.5cm,0.866cm>   
\setplotarea x from -10 to 12, y from -0.5  to 2.5         
\putrule from -2 1 to 6 1
\putrule from -4 0 to 4 0
\putrule from 1 1 to 1 0
\putrule from 6 1.2 to 6 0.8
\putrule from -2 1.2 to -2 0.8
\putrule from 4 0.2 to 4 -0.2
\putrule from -4 0.2 to -4 -0.2
\putrule from 5 2 to 10 2
\put{$_\bullet$}  at  -1 1
\put{$_{xO}$}[b]  at  -1 1.2
\put{$_\bullet$}  at  -3 0
\put{$_O$}[b]  at  -3 0.2
\put{$_e$}[r]  at  0.8 0.5
\put{$_T$}[l]  at  4.3 0
\put{$_{xT}$}[l]  at  6.3 1
\setlinear
\plot 4 1  5 2 /
\put{$_{\w}$}[l]  at  10.4  2
\arrow <6pt> [.3,.67] from  9.8 2 to  10 2 
\endpicture
}
\hfil
\caption{A boundary point $\w \in \Pi_x$.}
\end{figure}

\begin{proof} (a)
By (\ref{cov}), the element $\pi(x)p_{x^{-1}}\pi(x^{-1})$ is the characteristic function of the set
\begin{equation*}
\begin{split}
F_x& = \{ x\w \ ; \ \w\in\bD, x^{-1}T\cap [O,\w) \ne \emptyset \}  \\
& = \{ x\w \ ; \ \w\in\bD, T\cap [xO,x\w) \ne \emptyset \}  \\
& = \{ \w\in\bD \ ; \ T\cap [xO,\w) \ne \emptyset \}.
\end{split}
\end{equation*}
Now there exists a unique edge $e\in\D^1$ such that $o(e)\in T$ and $t(e)\in xT$. If $x\in S$ then $e\in\D^1_+$ and if $x\in S^{-1}$ then $e\in\overline{\D^1_+}$. Therefore
\begin{equation*}
\begin{split}
\bD-F_x& = \{ \w\in\bD \ ; \ T\cap [xO,\w) = \emptyset \}  \\
& = \{ \w\in\bD \ ; \ xT\cap [O,\w) \ne \emptyset \}  \\
& = \Pi_x ,
\end{split}
\end{equation*}
and the characteristic function of this set is $p_x$. See Figure \ref{fig1}.

The proof of (b) is an easy consequence of (\ref{cov}). 
\end{proof}

\refstepcounter{picture}

\begin{lemma}\label{L2}
The family of projections $P=\{p_g\ ; \ g\in \G-\{1\}\}$ generates $C(\bD)$ as a $C^*$-algebra.
\end{lemma}

\begin{figure}[htbp]
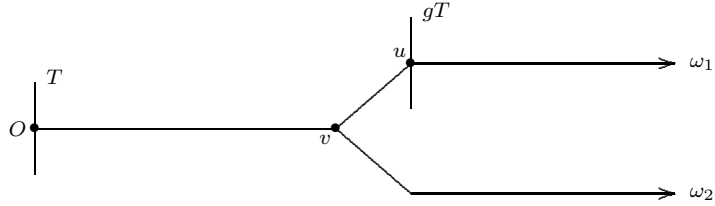
\label{fig2}
\hfil
\centerline{
\beginpicture
\setcoordinatesystem units <0.5cm,0.866cm>   
\setplotarea x from -10 to 10, y from -2  to 2         
\putrule from 0 0 to -8 0
\putrule from 2 1 to 9 1
\putrule from 2 -1 to 9 -1
\putrule from 2 1.7 to 2 0.3
\putrule from -8  0.7 to -8 -0.7
\setlinear
\plot 2 1  0 0  2 -1 /
\put{$_T$}[l]  at  -7.7 0.8
\put{$_{gT}$}[l]  at 2.3 1.8
\put{$_{\w_1}$}[l]  at  9.4  1
\put{$_{\w_2}$}[l]  at  9.4  -1
\put{$_\bullet$}  at  0 0
\put{$_\bullet$}  at  2 1
\put{$_\bullet$}  at  -8 0
\put{$_O$}[r]  at  -8.2 0
\put{$_{u}$}[b,r] at  1.9 1.1
\put{$_v$}[t,r] at  -0.1 -0.1
\arrow <6pt> [.3,.67] from  8.8 1 to  9 1
\arrow <6pt> [.3,.67] from  8.8 -1 to  9 -1
\endpicture
}
\hfil
\caption{Separation of boundary points.}
\end{figure}

\begin{proof}
We show that $P$ separates points of $\bD$. Let $\w_1, \w_2\in \bD$ with $\w_1\ne \w_2$. Let $[O,\w_1)\cap[O,\w_2)=[O,v]$, and choose $u\in[v,\w_1)$ such that $d(v,u)$ is greater than the diameter of $T$. See Figure \ref{fig2}.

Let $g\in\G$ be the unique element such that $u\in gT$. Then
$v\notin gT$ and so $gT\cap[O,\w_2)=\emptyset$.
Therefore $p_g(\w_1)=1$ and $p_g(\w_2)=0$.
\end{proof}

\refstepcounter{picture}

\begin{lemma}\label{L3}
The sets of the form $\Pi_x$, $x\in S\cup S^{-1}$, are pairwise disjoint and their union is $\bD$.
\end{lemma}

\begin{figure}[htbp]
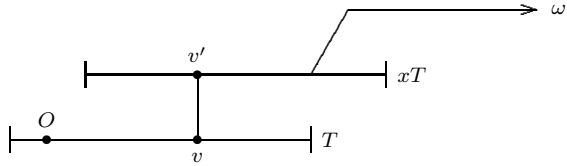
\label{fig1a}
\hfil
\centerline{
\beginpicture
\setcoordinatesystem units <0.5cm,0.866cm>   
\setplotarea x from -10 to 12, y from -0.5  to 2.5         
\putrule from -2 1 to 6 1
\putrule from -4 0 to 4 0
\putrule from 1 1 to 1 0
\putrule from 6 1.2 to 6 0.8
\putrule from -2 1.2 to -2 0.8
\putrule from 4 0.2 to 4 -0.2
\putrule from -4 0.2 to -4 -0.2
\putrule from 5 2 to 10 2
\put{$_\bullet$}  at  -3 0
\put{$_\bullet$}  at  1 0
\put{$_\bullet$}  at  1 1
\put{$_v$}[t]  at  1 -0.2
\put{$_{v^{\prime}}$}[b]  at 1 1.2 
\put{$_O$}[b]  at  -3 0.2
\put{$_T$}[l]  at  4.3 0
\put{$_{xT}$}[l]  at  6.3 1
\setlinear
\plot 4 1  5 2 /
\put{$_{\w}$}[l]  at  10.4  2
\arrow <6pt> [.3,.67] from  9.8 2 to  10 2 
\endpicture
}
\hfil
\caption{Definition of the set $\Pi_x$ containing $\w$.}
\end{figure}

\begin{proof}
Given $\w\in\bD$, let $v$ be the unique vertex of $\D$ such that 
$[O,\w)\cap T = [O,v]$. Let $v'$ be the vertex of $[O,\w)$ such that
$d(O,v')=d(O,v)+1$. Then let $x$ be the unique element of $S\cup S^{-1}$ such that $v'\in xT$. See Figure \ref{fig1a}. Then $\w\in \Pi_x$. The sets $\Pi_x$, $x\in S\cup S^{-1}$, are pairwise disjoint since the sets $xT$, $x\in S\cup S^{-1}$, are pairwise disjoint.
\end{proof}

\refstepcounter{picture}
\bigskip

For $x\in S\cup S^{-1}$ define a partial isometry
\[
s_x = \pi(x)(\id - p_{x^{-1}}) \in 
C(\bD)\rtimes\G.
\]
Then, by Lemma \ref{L1},
\[
s_xs_x^*= \pi(x)(\id - p_{x^{-1}})\pi(x^{-1})
= \id - \pi(x)p_{x^{-1}}\pi(x^{-1})
=p_x,
\]
and
\[
s_x^*s_x = \id - p_{x^{-1}}.
\]
Therefore the elements $s_x$ satisfy the relations
\begin{equation}\label{A}
s_x^*s_x = \displaystyle\sum_{\substack{y\in S\cup S^{-1} \\  y\ne x^{-1}}} s_ys_y^*.
\end{equation}
Also, it follows from Lemma \ref{L3} that
\begin{equation}\label{B}
\id =\displaystyle \sum_{x\in S\cup S^{-1}} p_x =
\displaystyle \sum_{x\in S\cup S^{-1}} s_xs_x^*.
\end{equation}

The relations (\ref{A}),(\ref{B}) are precisely the Cuntz-Krieger relations \cite{ck} corresponding to the \zomatrix \  $A$, with entries indexed by elements of $S\cup S^{-1}$, defined by
\begin{equation}\label{zo}
A(x,y) = 
\begin{cases}
    1&  \text{if $y\ne x^{-1}$},\\
    0&  \text{if $y= x^{-1}$.}
\end{cases}
\end{equation}
The matrix $A$ depends only on the rank of the free group $\G$.
Also $A$ is irreducible and not a permutation matrix. It follows that the $C^*$-subalgebra $\cA$ of $C(\bD)\rtimes\G$ generated by $\{s_x \ ; \ x\in S\cup S^{-1}\}$ is isomorphic to the simple Cuntz-Krieger algebra $\cO_A$ \cite{ck}.  It remains to show that 
$\cA$ is the whole of $C(\bD)\rtimes\G$.

\begin{lemma}\label{mainlemma}
Under the above hypotheses, $C(\bD)\rtimes\G = \cA$.
\end{lemma}

\begin{proof}
By the discussion above, it is enough to show that 
$$\cA\supseteq C(\bD)\rtimes\G.$$
First of all we show that $\cA\supseteq \pi(\G)$. 
It suffices to show that $\pi(x)\in \cA$ for each  $x\in S\cup S^{-1}$. Now
\[
s_{x^{-1}}^* = (\id - p_x)\pi(x)  =\pi(x) p_{x^{-1}},
\]
by Lemma \ref{L1}. Therefore
\begin{equation}\label{number}
s_x + s_{x^{-1}}^* =  \pi(x)(\id - p_{x^{-1}}) + \pi(x) p_{x^{-1}}
= \pi(x).
\end{equation}
It follows that $\pi(x)\in \cA$, as required.

Finally, we must show that $\cA\supseteq C(\bD)$. Since $s_xs_x^*=p_x$,
it is certainly true that $p_x\in \cA$ for all $x\in S\cup S^{-1}$.
It follows by induction from Lemma \ref{L1}(b), that $p_g\in \cA$ for all
$g\in \G$. 
Lemma \ref{L2} now implies that $\cA\supseteq C(\bD)$.
\end{proof}

\begin{example} \label{nonconjugate}
Consider the graphs $X$, $Y$ in Figure \ref{graphs}. Each of them has as universal covering space the 3-homogeneous tree $\D$.
Each has fundamental group the free group $\G$ on two generators.
Consequently, each gives rise to an action of $\G$ on $\D$. These two actions cannot be conjugate via an element of $\aut(\D)$ because their quotients are not isomorphic as graphs.

\refstepcounter{picture}
\begin{figure}[htbp]
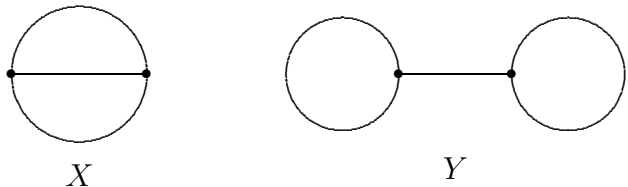
\label{graphs}
\hfil
\centerline{
\beginpicture 
\setcoordinatesystem units <0.3cm, 0.3cm>    
\setplotarea  x from -4 to 4,  y from -4.5 to 4
\put{$_\bullet$}  at -3 0
\put{$_\bullet$}  at 3 0
\putrule from -3 0 to 3 0
\setlinear
\circulararc 360 degrees from 3 0 center at 0 0
\put{$X$}  at  0 -4.5 
\setcoordinatesystem units <0.5cm, 0.5cm>  point at -10 0  
\setplotarea  x from -4.5 to 4.5,  y from -2.5 to 2
\putrule from -1.5 0 to 1.5 0
\put{$_\bullet$}  at  -1.5 0
\put{$_\bullet$}  at  1.5 0
\put{$Y$}  at  0  -2.5
\circulararc 360 degrees from 3 -1.5 center at 3 0 
\circulararc 360 degrees from -3 1.5 center at -3 0
\endpicture
}
\hfil
\caption{The graphs $X$, $Y$.}
\end{figure}

Copies of the free group on two generators, acting with these actions on the corresponding Bruhat-Tits tree $\D$, can be found inside $\PGL(2,\QQ_2)$, and also inside $\PGL(2,\FF)$ for any local field $\FF$ with residue field of order $2$. To see this, note that by \cite[Appendix, Proposition 5.5]{ftn} $\PGL(2,\QQ_2)$ contains cocompact lattices 
$\G_1$, $\G_2$ which act freely and transitively on the vertex set $\D^0$ with
\begin{equation*}
\begin{split}
\G_1& = (\ZZ/2\ZZ)*(\ZZ/2\ZZ)*(\ZZ/2\ZZ)
=\langle a, b, c \, |\, a^2=b^2=c^2=1 \rangle \\
\G_2& = \ZZ * (\ZZ/2\ZZ) 
=\langle x, d \, |\, d^2=1 \rangle.
\end{split}
\end{equation*}
The subgroups $\G_X=\langle ab, ac \rangle$ and $\G_Y=\langle x, dxd \rangle$
are both isomorphic to the free group on two generators. Moreover 
$\G_X\backslash \D = X$ and $\G_Y\backslash \D = Y$. 
\end{example}

\bigskip

\section{K-theory}\label{Ksection}

Using the results of \cite{c}, it is now easy to determine the K-theory of $\cA(\G)$.
For each $x\in S\cup S^{-1}$, the element $p_x$ is a projection in $\cA(\G)$ and therefore defines an equivalence class $[p_x]$ in $K_0(\cA(\G))$. 
It is shown in \cite{c} that the classes $[p_x]$ generate $K_0(\cA(\G))$.
Indeed, let $L$ denote the abelian group with generating set $S\cup S^{-1}$ and relations

\begin{equation}\label{rel1}
x = \displaystyle\sum_{\substack{y\in S\cup S^{-1} \\  y\ne x^{-1}}}y  \qquad \text{for} \ x\in S\cup S^{-1}. 
\end{equation}

The map $x\mapsto [p_x]$ extends to an isomorphism $\theta$ from $L$ onto $K_0(\cA(\G))$ \cite{c}. Moreover $\theta(\e)=[\id]$, where 
$\e= \displaystyle\sum_{x\in S\cup S^{-1}}x$.  Now it follows from (\ref{rel1}) that,
for each $x\in S$,
\begin{equation}\label{rel2}
\e = x + x^{-1}. 
\end{equation}
Also
\begin{equation*}
\e = \displaystyle\sum_{x\in S}(x + x^{-1})=\displaystyle\sum_{x\in S}\e =\g\e. 
\end{equation*}
Thus
\begin{equation}\label{rel3}
(\g-1)\e = 0. 
\end{equation}

The group $L$ is therefore generated by $S\cup\{\e\}$, and the relation (\ref{rel3})
is satisfied.

On the other hand, starting with an abstract abelian group with generating set  $S\cup\{\e\}$ and the relations (\ref{rel2}), one can make the formal definition $x^{-1}=\e-x$, for each $x\in S$, and recover the relations (\ref{rel1}) via
\begin{equation*}
\displaystyle\sum_{x\in S}(x + x^{-1})= \g \e = \e = x+ x^{-1} \qquad \text{for} \ x\in S. 
\end{equation*}
This discussion proves

\begin{lemma}\label{Ktheory}
$K_0(\cA(\G))\cong \ZZ^\g \oplus \ZZ/(\g -1)\ZZ$ via an isomorphism which sends $[\id]$ to the generator of $\ZZ/(\g -1)\ZZ$.
\end{lemma}

It is known that the $C^*$-algebra $\cA(\G)$ is purely infinite, simple, unital and nuclear \cite{ck,c,c'}. The classification theorem of \cite{k} therefore shows that $\cA(\G)$ is determined by its K-theory.

\begin{remark}\label{forward}
If $\G$ is a torsion free cocompact lattice in $\PSL (2,\RR)$, so that $\G$ is the fundamental group of a Riemann surface of genus $g$, then it is known, \cite[Proposition 2.9]{ad}, \cite{hn}, that $\cA(\G)=C(\PP_1(\RR))\rtimes \G$ is the unique p.i.s.u.n. $C^*$-algebra whose K-theory is specified by
\begin{align*}
(K_0(\cA(\G)), [\id]) &=(\ZZ^{2g+1}\oplus\ZZ/(2g-2)\ZZ, (0,0,\dots ,0,1)), \\
K_1(\cA(\G)) &= \ZZ^{2g+1}.
\end{align*}

The proof of this result in \cite{ad} makes use of the Thom Isomorphism Theorem of A. Connes (which has no $p$-adic analogue) to identify  $K_*(\cA(\G))$ with the topological K-theory  $K^*(\G\backslash \PSL (2,\RR))$.
It follows from the classification theorem of \cite{k,ph} that $\cA(\G)$ is a Cuntz-Krieger algebra. However there is no apparent dynamical reason for this fact. In contrast, the Cuntz-Krieger algebras of the present article appear naturally and explicitly.
\end{remark}

\bigskip


\section{The measure theoretic result}\label{measure section}

The purpose of this section is to prove Theorem \ref{main2} of the Introduction.
From now on $\D$ is a homogeneous tree of degree $q+1$, where $q\ge 1$, and $\G$ is a free uniform lattice in $\aut(\D)$. A similar Theorem could be stated for non-homogeneous trees, and proved by the same methods.  The boundary $\bD$ is endowed with  a natural Borel measure. In contrast to the topological result, measure theoretic rigidity for the boundary action fails: the von Neumann algebra $\Li(\bD)\rtimes \G$ depends on the tree $\D$ and on the action of $\G$.  Before proceeding with the proof here are some examples.

\begin{example}\label{measureexample1}
Let $\G$ be the free group on two generators. Then $\G$ is the fundamental group of each of the graphs $X$, $Y$ of Figure \ref{graphs}. The 3-homogeneous tree $\D_3$ is the universal covering of both these graphs and there are two corresponding (free, cocompact) actions of $\G$ on $\D_3$.
It follows from Theorem \ref{main2} that the von Neumann algebra $\Li(\bD_3)\rtimes \G$ is the hyperfinite factor of type
${{\text{\rm{III}}}_{1/4}}$ in the first case, since $X$ is bipartite, and type
${{\text{\rm{III}}}_{1/2}}$ in the second case, since $Y$ is not bipartite.

The group $\G$ is also the fundamental group of a bouquet of two circles and the corresponding action of $\G$ on the 4-homogeneous tree $\D_4$ produces the hyperfinite factor of type ${{\text{\rm{III}}}_{1/3}}$.  These three actions are the only free and cocompact actions of the free group on two generators on a tree $\D$ with no vertices of degree $\le 2$. 
\end{example}

\begin{remark}
For each $\g\ge 2$, it is easy to construct bipartite and non-bipartite 3-homogeneous graphs with fundamental group the free group on $\g$ generators. The corresponding boundary actions are of types ${{\text{\rm{III}}}_{1/4}}$ and ${{\text{\rm{III}}}_{1/2}}$
respectively. 
\end{remark}

We now proceed with the proof of Theorem \ref{main2}. As before, fix a vertex $O\in \D$. If $u, v$ are vertices in $\D^0$, let $[u,v]$ be the directed geodesic path between them, with origin $u$. The graph distance $d(u,v)$ between $u$ and $v$ is the length of $[u,v]$, where each edge is assigned unit length. If $v\in \D^0$ let $\Om_v$ be the clopen set consisting of all $\w \in\bD$ such that $v\in [O,\w)$.

\refstepcounter{picture}
\begin{figure}[htbp]
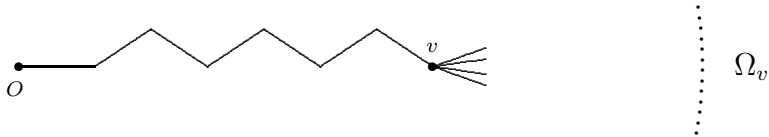
\label{t(w)}
\hfil
\centerline{
\beginpicture
\setcoordinatesystem units <0.5cm, 0.5cm>  
\setplotarea  x from -8 to 10.5,  y from -2 to 3
\put{$_{O}$}   [t] at -8.1 -0.4
\put{$_{\bullet}$}    at -8 0
\put{$_{\bullet}$}    at 3 0
\put{$_{v}$}   [b] at 3  0.4
\put{$\Omega_v$}   [l] at 11 0
\putrule from -8 0 to -6 0
\setlinear
\plot   -6 0  -4.5  1   -3  0  -1.5 1   0  0  1.5 1   3  0  4.4 0.5 /
\plot    3  0  4.4 0.2 /
\plot    3  0  4.4 -0.2 /
\plot    3  0  4.4 -0.5 /
\setplotsymbol({$\cdot$}) \plotsymbolspacing=4pt
\circulararc 20 degrees from 10 -1.8 center at 0 0 
\endpicture
}
\hfil
\caption{A subset $\Om_v$ of the boundary.}
\end{figure}

There is a natural Borel measure $\mu$ on $\bD$ defined by $\mu (\Om_v)=q^{(1-n)}$, where $n=d(O,v)$. The consistency of this definition is easily established, using the fact that there are precisely $q$ vertices $w$ adjacent to $v$ and  not lying on the path $[O,v]$. The set $\Om_v$ is the disjoint union of the corresponding sets $\Om_w$, each of which has measure $q^{-n}$. Note that the normalization of this measure is different from that in \cite{ftn}. This is immaterial for the result, but makes the formulae simpler. The measure $\mu$ clearly depends on the choice of the vertex $O$ in $\D$, but its measure class does not.

\begin{lemma}\label{measure-theoretic freeness}
The action of $\G$ on $\bD$ is measure-theoretically free, i.e.
\[
\mu\left(\{ \w\in\bD : g\w= \w\}\right) = 0
\]
for all elements $g\in\G -\{ e\}$.
\end{lemma}

\begin{proof}
Let $g\in\G-\{ e\}$. Since the action of $\G$ on $\D$ is free,
$g$ is hyperbolic; that is $g$ fixes no point of $\D$.
It follows that the set
$\{\w\in\bD : g\w= \w\}$
contains exactly two elements and so certainly has measure zero.
\end{proof}

It is well known (and it is an easy consequence of Lemma \ref{ergodicH} below) that the action of $\G$ on $\bD$ is also ergodic. Therefore the von Neumann algebra $\Li(\bD)\rtimes \G$ is a factor. A convenient reference for this fact and for the classification of von Neumann algebras is \cite{su}.  Most of this section will be devoted to establishing that this factor is of type $\tl$, for an appropriate value of $\lambda$. This will be done by determining the ratio set of W.~Krieger.

\begin{definition}
Let $G$ be a countable group of automorphisms of a measure space $(\Om,\mu)$.
Define the {\bf ratio set} $r(G)$ to be the subset of
$[0,\infty)$ such that if $\lambda\geq 0$ then $\lambda\in r(G)$ if and
only if for
every $\epsilon>0$ and measurable set $A$ with $\mu(A)>0$, there exists $g \in G$ and a measurable set $B$ such that $\mu(B)>0$, $B\cup gB\subseteq A$ and
\[
\left| \frac{d\mu{\scriptstyle\circ}g}{d\mu}(\w) -\lambda \right| <\epsilon
\]
for all $\w\in B$.
\end{definition}

\begin{remark}
The ratio set $r(\G)$ depends only on the quasi-equivalence class of the measure
$\mu$. If the action of $\G$ is ergodic then $r(\G)-\{0\}$ is a subgroup of the multiplicative group of positive real numbers \cite[\S I-3, Lemma 14]{ho}.
\end{remark}

In order to compute $r(\G)$, for the action of $\G$ on $\bD$, the first step is to find the possible values of the Radon-Nikodym derivatives 
$\frac{d\mu{\scriptstyle\circ}g}{d\mu}(\w)$, for $g\in \G$ and $\w\in \bD$.

Fix $g\in \G$ and $\w\in \bD$.
Choose an open set of the form  $\Om_v$ with $v\in [O,\w)$ and  $d(O,v)>d(O,gO)$.  Such sets $\Om_v$ form a neighbourhood base of $\w$.
Then $v\notin [O,gO]$ (Figure \ref{basicnbd}), and $g^{-1}\Om_v=\Om_{g^{-1}v}$.
Since $d(O,g^{-1}v)=d(gO,v)$, we have  $\mu(g^{-1}\Om_v)=q^{-d(gO,v)}$.
\refstepcounter{picture}
\begin{figure}[htbp]\label{basicnbd}
\hfil
\centerline{
\beginpicture
\setcoordinatesystem units <0.5cm,0.866cm>   
\setplotarea x from -10 to 10, y from -0.5  to 1         
\putrule from -8 0 to 9 0
\setlinear
\plot  -6 0  -4 1 /
\plot    2  0  3.4 0.5 /
\plot    2  0  3.4 0.2 /
\plot    2  0  3.4 -0.2 /
\plot    2  0  3.4 -0.5 /
\put{$_{\w}$}[l]  at  9.4  0
\put{$_\bullet$}  at  -4  1
\put{$_{gO}$}[l]  at  -3.6  1
\put{$_\bullet$}  at  -8  0
\put{$_O$}[t]  at  -8  -0.2
\put{$_\bullet$}  at  2   0
\put{$_{v}$}[t]  at  2   -0.3
\arrow <6pt> [.3,.67] from  8.8 0 to  9 0
\endpicture
}
\hfil
\caption{}
\end{figure}

It follows that
\begin{equation}\label{rnn}
\frac{d\mu{\scriptstyle\circ}g}{d\mu}(\w)=\frac {\mu(g^{-1}\Om_v)}{\mu(\Om_v)}
=\frac {q^{-d(gO,v)}}{q^{-d(O,v)}}=q^{\delta(g,\w)}
\end{equation}
where $\delta(g,\w)=d(O, v)-d(gO, v)$. It is clear that $\delta(g,\w)$ depends only on $g$ and  $\w$, not on the choice of $v$. In the language of \cite[Chapter 8]{gh}, $\delta(g,\w)$ is the Busemann function $\beta_{\w}(O,gO)$ relating the horocycles centered at $\w$ containing $O$, $gO$ respectively.  For a fixed vertex $v$ with $d(O,v)>d(O,gO)$, the formula (\ref{rnn}) remains true for all $\w\in\Om_v$.

We have therefore proved

\begin{lemma}\label{first}
The values of the Radon-Nikodym derivatives 
$\frac{d\mu{\scriptstyle\circ}g}{d\mu}(\w)$, for $g\in \G$ and $\w\in\bD$, are given by
\[
\frac{d\mu{\scriptstyle\circ}g}{d\mu}(\w)=q^{\delta(g,\w)}
\]
Moreover, for each $g\in \G$, each of these values is attained on a nonempty open subset of $\bD$.
\end{lemma}

These considerations show  that
\begin{equation}\label{possibleRN}
r(\G)\subseteq \{ q^{\delta(g,\w)} \, ;\, g\in \G, \w\in\bD \}\cup\{0\}.
\end{equation}

Since the action of $\G$ is ergodic, $r(\G)-\{0\}$ is a multiplicative group of positive real numbers \cite[Lemma 14]{ho}.
What must be done now is to show that the inclusion in (\ref{possibleRN}) is in fact an equality.
Clearly $r(\G)\ne [0,\infty)$. Therefore if we can show that $r(\G)$ contains a number in the open interval $(0,1)$ then, by \cite[Lemma 15]{ho}, it must equal $\{ \lambda^n \, ;\,  n\in \ZZ \}\cup\{0\}$, for some $\lambda \in (0,1)$. By definition, this will show that the action of $\G$, and hence the associated von Neumann algebra $\Li(\bD)\rtimes \G$), is of type $\tl$.

\smallskip

\begin{figure}[htbp]
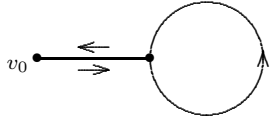
\label{cycle with tail}
\hfil
\centerline{
\beginpicture 
\setcoordinatesystem units <0.5cm, 0.5cm>  
\setplotarea  x from -3 to 3,  y from -1.5 to 1.5
\putrule from -1.5 0 to 1.5 0
\put{$_{v_0}$}[r]  at  -1.7 -0.2
\put{$_\bullet$}  at  -1.5 0
\put{$_\bullet$}  at  1.5 0
\circulararc 360 degrees from 3 -1.5 center at 3 0 
\arrow <6pt> [.3,.67] from  -0.4 -0.3 to  0.4 -0.3
\arrow <6pt> [.3,.67] from   0.4 0.3 to  -0.4 0.3
\arrow <6pt> [.3,.67] from  4.5  0 to  4.5 0.2
\endpicture
}
\hfil
\caption{A proper cycle with tail, based at $v_0$.}
\end{figure}

Before proceeding, it is useful to interpret the situation in terms of the quotient graph $X=\G\backslash \D$.
In a connected graph $X$ a {\em proper} path is a path which has no backtracking. That is, no edge $[a,b]$ in the path is immediately followed by its inverse $[b,a]$. 
A {\em cycle} is a closed path, which is said to be based at its initial vertex (= final vertex).
Note that a proper cycle can have a tail beginning at its base vertex, but that it can have no other tail (Figure \ref{cycle with tail}).  Every proper cycle determines a unique tail-less cycle which is obtained by removing the tail.  A circuit is a cycle which does not pass more than once through any vertex. There is clearly an upper bound for the possible length of a circuit in $X$. 

\refstepcounter{picture}

If $g\in \G$, then the geodesic path $[O, gO]$ in $\D$ projects to a proper cycle $C$ in the quotient graph $X=\G\backslash\D$ based at $v_0=\G O$. Moreover $d(O, gO)$ is equal to the length $\ell (C)$ of that cycle.

Conversely if $C$ is a proper cycle based at $v_0$ in the graph $X$ then the homotopy class of $C$ is an element $g$ of the fundamental group $\G$ of $X$. The cycle $C$ lifts to a unique proper path in $\D$ with initial vertex $O$, namely  $[O, gO]$, and
$\ell (C)=d(O, gO)$.

\medskip

In order to prove equality in (\ref{possibleRN}) we need the auxiliary concept of the {\em full group}.

\begin{definition}
Given a group $\G$ acting on a measure space $(\Om,\mu)$, we define the {\em full
group}, $[\G]$, of $\G$ by
\[
[\G]=\left\{ T\in\aut(\Om) : T \w\in \G\w \text{ for almost every
}\w\in\Om\right\}.
\]
\end{definition}

\begin{remark}\label{full group determines ratio set}
The ratio set $r(\G)$ of a countable subgroup $\G$ of $\aut(\Om)$ depends only on the full
group in the sense that $r(\G_1)=r(\G_2)$ whenever $[\G_1]=[\G_2]$.
\end{remark}

The basis for the proof of equality in (\ref{possibleRN}) is the following well known result. It is stated without proof in \cite[I.3]{ho}.

\begin{lemma}\label{secondratio}
Let $\G$ be a countable group acting ergodically on a measure space $\Om$. Suppose that
the full group $[\G]$ contains an ergodic measure preserving subgroup $H$.

If $r \in (0,\infty)$,  $g\in \G$ and the set  
$D=\{\w\in \Om \, ;\, \frac{d\mu{\scriptstyle\circ}g}{d\mu}(\w) = r\}$ has positive measure, then $r\in r(\G)$.
\end{lemma}

\begin{proof}
Let $A$ be a measurable subset of $\Om$ with $\mu(A)>0$.
By the ergodicity of $H$,
there exist $h_1,h_2\in H$ such that the set
$B=\{\w\in A: h_1\w\in D \text{ and } h_2gh_1\w\in A\}$
has positive measure.

Let $\G^{\prime}$ denote the group generated by $h_1, h_2$ and $\G$.
By Remark \ref{full group determines ratio set}, $r(\G^{\prime}) =r(\G)$.

Let $t=h_2gh_1\in \G^{\prime}$. By construction, $B\cup tB\subseteq A$.
Moreover, since $H$ is measure-preserving,
\[
\frac{d\mu{\scriptstyle\circ}t}{d\mu}(\w)=
\frac{d\mu{\scriptstyle\circ}g}{d\mu}(h_1\w)= r
 \text{ for all } \w\in B,
\]
 since $h_1\w \in D$. This proves that $r\in r(\G^{\prime}) = r(\G)$, as required.
\end{proof}

In view of Lemma \ref{first}, all that is now needed in the present setup is
the construction of a subgroup $H$. This will require the following result.

\begin{lemma}\label{mpergodicity}
Let $H$ be subgroup of $\aut(\D)$. Suppose that the induced action of $H$ on $\bD$ is measure preserving and that, for each positive integer $n$, $H$ acts transitively on the
collection of sets 
\[
\left\{\Om_v :  v\in \D ,\, d(O,v)=n \right\}.
\] 
Then $H$ acts ergodically on $\bD$.
\end{lemma}
\begin{proof}
Suppose that $S_0\subseteq\bD$ is a Borel set which is invariant under $H$ and
such that $\mu(S_0)>0$. We show that this implies
$\mu(\bD- S_0)=0$, thereby establishing the ergodicity of the action.

Define a new measure $\lambda$ on $\bD$ by $\lambda(S)=\mu(S\cap S_0)$, for each Borel
set $S\subseteq\bD$. Now, for each $k\in H$,
\begin{eqnarray*}
\lambda(kS)&=&\mu(kS\cap S_0) = \mu(S\cap k^{-1}S_0) \\
&\leq & \mu(S\cap S_0) + \mu(S\cap (k^{-1}S_0- S_0)) \\
&=& \mu(S\cap S_0) \\
&=& \lambda(S),
\end{eqnarray*}
and therefore $\lambda$ is $H$-invariant.

Fix a positive integer $n$. The transitivity hypothesis on the action of $H$ implies that
\[
\lambda(\Om_v)=\lambda(\Om_w)
\]
whenever $v,w\in \D$, $d(O,v)=d(O,w)=n$. Since $\bD$ is the union of $q^{(n-1)}(q+1)$ disjoint sets
$\{\Om_v$; $d(O,v)=n\}$, each of which has equal measure with respect to $\lambda$, we deduce that, if $d(O,v)=n$,
\[
\lambda(\Om_v)=\frac{\lambda(\bD)}{q^{(n-1)}(q+1)}
=\frac{\mu(S_0)}{q^{(n-1)}(q+1)}.
\]
Thus $\lambda(\Om_v)=c\mu(\Om_v)$ for every $v\in \D$, 
where $c=\frac{\mu(S_0)}{(q+1)} >0$.
Since the sets $\Om_v$, $v\in \D$, generate the Borel $\s$-algebra, we deduce that
$\lambda(S)=c\mu(S)$ for each Borel set $S$. Therefore
\[
\mu(\bD- S_0) = c^{-1}\lambda(\bD- S_0) 
= c^{-1}\mu((\bD- S_0) \cap S_0) = 0 ,
\]
thus proving ergodicity.
\end{proof}

\medskip

It is now convenient to introduce some new terminology.

\begin{definition}\label{fL}
Let $X$ be a finite connected graph.
Let $v_0$ be a vertex of $X$ and let $K \ge 0$.
Say that $(X, v_0)$ has property $\fL(K)$ if for any two proper paths $P_1, P_2$ having the same length $n$ and the same initial vertex $v_0$, there exists $k\ge 0$, with $k\le K$, and proper cycles $C_1, C_2$ based at $v_0$ such that 
\begin{itemize}
\item[(a)] The initial segment of $C_i$ is $P_i$, $i=1,2$;
\item[(b)] the cycles $C_i$ have the same length $n+k$, $i=1,2$.
\end{itemize}
\end{definition}
\smallskip
Property $\fL(K)$ says that any two proper paths of the same length starting at $v_0$ can be completed to proper cycles of the same length, with a uniform bound on how much must be added to each path.
\medskip

\begin {lemma}\label{bcondition}
Let $X$ be a finite connected graph whose vertices all have degree at least three and let $v_0$ be a vertex of $X$.
Then $(X, v_0)$ has property $\fL(K)$ for some $K\ge 0$.
\end{lemma}

The proof of this technical result is deferred to Section \ref{technicallemma}.
We can now prove that the action of $\G$ on $\bD$ satisfies the hypotheses of Lemma \ref{secondratio}.

\begin {lemma}\label{ergodicH}
Let $\D$ be a homogeneous tree of degree $q+1$, where $q\ge 1$ and let $\G$ be a free uniform lattice in $\aut(\D)$.
Then, relative to the action of $\G$ on $\bD$, the full group $[\G]$  contains an ergodic measure preserving subgroup $H$.
\end {lemma}

\begin{proof}
By Lemma \ref{mpergodicity}, it suffices to prove the following assertion for any $u,v\in \D^0$, with $d(O,u)=d(O,v)=n$.
\begin{itemize}
\item[$(\star)$] There exists a measure preserving automorphism $\phi\in[\G]$ such that $\phi$ is almost everywhere a bijection from $\Om_u$ onto $\Om_v$.
\end{itemize}

The geodesic paths $[O,u], [O,v]$ in $\D$ project to proper paths  $P_u, P_v$ in $X$ with initial vertex $v_0=\G O$ and length $n$.
By hypothesis, the graph $X$ has property $\fL(K)$ for some constant $K\ge 0$, relative to $v_0$.  Therefore 
there exists an integer $k\le K$ and proper cycles $C_u, C_v$ based at $v_0$ which have initial segments $P_u, P_v$ respectively
and $\ell(P_u)=\ell(P_v)=n+k$. 

The cycles $C_u, C_v$ lift to unique geodesic paths $[O,u^*], [O,v^*]$ in $\D$ with initial segments $[O,u]$, $[O,v]$ respectively and $\G u^*=\G v^* = \G O = v_0$.
Since the vertices of $X$ all have degree at least three, we can choose an edge $e$ with $o(e)=v_0$ in $X$ such that $e$ meets the terminal edges of $C_u$ and $C_v$ only at $v_0$. There are unique vertices $u_1, v_1 \in \D^1$ such that  $e=\G [u^*,u_1] = \G [v^*,v_1]$. Therefore there exists an element $g\in \G$ such that 
$g[u^*,u_1] = [v^*,v_1]$. The restriction of the action of $g$ to 
$\Om_{u_1}=\{\w\in \bD ;\, u_1\in [u^*,\w)\}$ defines a measure preserving bijection from $\Om_{u_1}$ onto $\Om_{v_1}$.
Define $\phi(\w)= g(\w)$ for $\w\in \Om_{u_1}$. 

\refstepcounter{picture}
\begin{figure}[htbp]
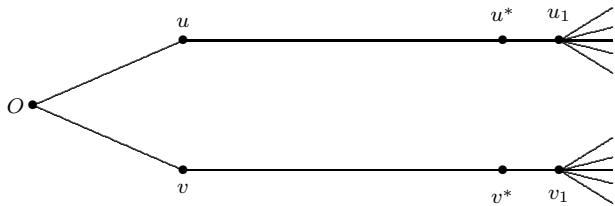
\label{partialdef}
\hfil
\centerline{
\beginpicture
\setcoordinatesystem units <0.5cm,0.866cm>   
\setplotarea x from -12 to 4, y from -3 to 1         
\putrule from -8 0 to 3.4 0
\putrule from -8 -2 to 3.4 -2
\setlinear
\plot  -8 -2  -12 -1  -8 0 /
\plot    2  0  3.4 0.5 /
\plot    2  0  3.4 0.2 /
\plot    2  0  3.4 -0.2 /
\plot    2  0  3.4 -0.5 /
\plot    2  -2  3.4 -1.5 /
\plot    2  -2  3.4 -1.8 /
\plot    2  -2  3.4 -2.2 /
\plot    2  -2  3.4 -2.5 /
\put{$_\bullet$}  at  -12 -1
\put{$_{O}$}[r]  at  -12.2 -1
\put{$_\bullet$}  at  -8  0
\put{$_u$}[b]  at  -8  0.2
\put{$_\bullet$}  at  0.5   0
\put{$_{u^*}$}[b]  at  0.5   0.3
\put{$_\bullet$}  at  2   0
\put{$_{u_1}$}[b]  at  2   0.3
\put{$_\bullet$}  at  -8  -2
\put{$_v$}[t]  at  -8  -2.2
\put{$_\bullet$}  at  0.5   -2
\put{$_{v^*}$}[t]  at  0.5   -2.3
\put{$_\bullet$}  at  2   -2
\put{$_{v_1}$}[t]  at  2   -2.3
\endpicture
}
\hfil
\caption{Definition of $\phi$ on $\Om_{u_1}$.}
\end{figure}

The set $\Om_u$ is a disjoint union of $q^{k+1}$ sets of the form $\Om_w$ where $d(O,w)=n+k+1$. Each such set therefore has measure $\mu(\Om_w)=q^{-k-1}\mu(\Om_u)$.
The map $\phi$ has been defined only on the set $\Om_w$ with $w=u_1$.
Therefore $\phi$ has not yet been defined on a proportion $(1-q^{-k-1})$ of the set $\Om_u$.  
Since $k\le K$, the measure of the subset of $\Om_u$ for which $\phi$ has not yet been defined is at most $(1-q^{-K-1})\mu(\Om_u)$.

\refstepcounter{picture}
\begin{figure}[htbp]
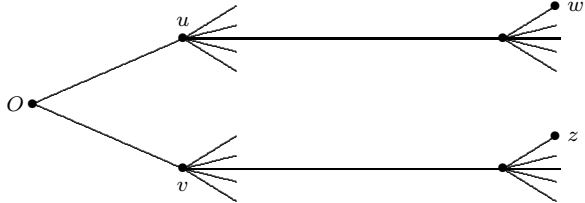
\label{notdef}
\hfil
\centerline{
\beginpicture
\setcoordinatesystem units <0.5cm,0.866cm>   
\setplotarea x from -12 to 3, y from -3 to 1         
\putrule from -8 0 to 2 0
\putrule from -8 -2 to 2 -2
\setlinear
\plot  -8 -2  -12 -1  -8 0 /
\plot    0.5  0  1.9 0.5 /
\plot    0.5  0  1.9 0.2 /
\plot    0.5  0  1.9 -0.2 /
\plot    0.5  0  1.9 -0.5 /
\plot    0.5  -2  1.9 -1.5 /
\plot    0.5  -2  1.9 -1.8 /
\plot    0.5  -2  1.9 -2.2 /
\plot    0.5  -2  1.9 -2.5 /
\plot    -8  0  -6.6 0.5 /
\plot    -8  0  -6.6 0.2 /
\plot    -8  0  -6.6 -0.2 /
\plot    -8  0  -6.6 -0.5 /
\plot    -8  -2  -6.6 -1.5 /
\plot    -8  -2  -6.6 -1.8 /
\plot    -8  -2  -6.6 -2.2 /
\plot    -8  -2  -6.6 -2.5 /
\put{$_\bullet$}  at  1.9 0.5
\put{$_w$}[l]  at  2.2 0.5
\put{$_\bullet$}  at  1.9 -1.5
\put{$_z$}[l]  at  2.2 -1.5
\put{$_\bullet$}  at  -12 -1
\put{$_{O}$}[r]  at  -12.2 -1
\put{$_\bullet$}  at  -8  0
\put{$_u$}[b]  at  -8  0.2
\put{$_\bullet$}  at  0.5   0
\put{$_\bullet$}  at  -8  -2
\put{$_v$}[t]  at  -8  -2.2
\put{$_\bullet$}  at  0.5   -2
\endpicture
}
\hfil
\caption{Second step in the definition of $\phi$.}
\end{figure}

Now repeat the process on each of the $q^{k+1}$ sets $\Om_w\subset \Om_u$, with $d(O,w)=n+k+1$, on which $\phi$ has not yet been defined. In the preceding argument,  replace $\Om_u$ by $\Om_w$ and $\Om_v$ by an appropriate subset $\Om_z$ of $\Om_v$ disjoint from $\Om_{v_1}$. Note that there is a large amount of arbitrariness in the choice of which $\Om_z$ is to be paired with a particular $\Om_w$.  
In each such $\Om_w$, $\phi$ is then defined on a subset whose complement in $\Om_w$ has measure at most $(1-q^{-K-1})\mu(\Om_w)$.

Thus after two steps, $\phi$ has been defined except on a set of measure at most $(1-q^{-K-1})^2\mu(\Om_u)$.  Continue in this way. After  $j$ steps, $\phi$ has been defined except on a set of measure at most $(1-q^{-K-1})^j\mu(\Om_u)$.

Since $(1-q^{-K-1})^j\to 0$ as $j\to\infty$, the measure preserving map $\phi$ is defined almost everywhere on $\Om_u$, with $\phi(\w)\in \G\w$ for almost all $\w\in \Om_u$.  Finally define $\phi$ to be the inverse of the map already constructed on $\Om_v$ and the identity map on $\bD-(\Om_u\cup\Om_v)$.  The proof of $(\star)$ is complete.
\end{proof}

It follows from Lemmas \ref{first}, \ref{secondratio}, and \ref{ergodicH} that we have equality in (\ref {possibleRN}). That is
\begin{equation}\label{exactRN}
r(\G)= \{ q^{\delta(g,\w)} \, ;\, g\in\G, \w\in \bD \}\cup\{0\}.
\end{equation}

The final step is to identify this set more precisely.
Recall that there is a canonical bipartition of the vertex set of $\D$, such that two vertices have the same type if and only if the distance between them is even. The graph $X=\G\backslash\D$ is bipartite if and only if the action of $\G$ is type preserving.
Recall also that $\delta(g,\w)=d(O, v)-d(gO, v)$, for any vertex $v\in [O,\w)\cap[gO,\w)$.

\begin{lemma}\label{gcd} Let $\D$ be a locally finite tree whose vertices all have degree at least three. Let $\G$ be a free uniform lattice in $\aut(\D)$ and let $X=\G\backslash\D$. Then
\begin{equation*}
\{ \delta(g,\w) \, ;\, g\in \G,\w\in\bD \}= 
\begin{cases}
2\ZZ &  \text{if $X$ is bipartite}, \\    
\ZZ &  \text{otherwise}.
\end{cases}
\end{equation*}
\end{lemma}

\begin{proof} Suppose first of all that $X$ is not bipartite. Then $X$ contains a circuit of odd length. Connecting this circuit to $v_0$ by a minimal path and going around the circuit an appropriate number of times shows that $X$ contains proper cycles based at $v_0$ of arbitrarily large even and odd lengths.

It follows that we may choose $g\in \G$ such that $d(O,gO)=2n$, for arbitrarily large $n$. If $k\in \ZZ$,with $-n\le k\le n$, let $a \in [O,gO]$ with $d(O,a)=n+k$, $d(gO,a)=n-k$. Choose $\w\in\bD$ with $[O,\w)\cap[gO,\w)=[a,\w)$. This is possible since the vertex $a$ has degree at least three (Figure \ref{delta}). Then 
$\delta(g,\w)= n+k -(n-k) = 2k$.

We may also choose $g\in \G$ such that $d(O,gO)=2n+1$ for arbitrarily large $n$.
 If $k\in \ZZ$,with $-n\le k\le n$, choose $a \in [O,gO]$ with $d(O,a)=n+k$, $d(gO,a)=n+1-k$. Choose $\w\in\bD$ with $[O,\w)\cap[gO,\w)=[a,\w)$. Then $\delta(g,\w)=2k-1$.

It follows that the range of the function $\delta(g,\w)$ is $\ZZ$.

\refstepcounter{picture}
\begin{figure}[htbp]\label{delta}
\hfil
\centerline{
\beginpicture
\setcoordinatesystem units <0.5cm,0.866cm>   
\setplotarea x from -3 to 10, y from -2  to 1         
\putrule from 2 0 to 9 0
\setlinear
\plot  -1 1  2 0 /
\plot  -4 -2  2 0 /
\put{$_{\w}$}[l]  at  9.4  0
\put{$_\bullet$}  at  -1 1
\put{$_O$}[r]  at  -1.2 1
\put{$_{gO}$}[r]  at  -4.2 -2
\put{$_\bullet$}  at  -4 -2
\put{$_{a}$}[t]  at  2   -0.3
\put{$_\bullet$}  at  2 0
\arrow <6pt> [.3,.67] from  8.8 0 to  9 0
\endpicture
}
\hfil
\caption{}
\end{figure}

Now suppose that $X$ is bipartite. Then the graph $X$ contains proper cycles of arbitrarily large even length only.  The preceding argument shows that the range of the function $\delta(g,\w)$ is $2\ZZ$.
\end{proof}

\smallskip

Assume the hypotheses of Theorem \ref{main2}. Then by (\ref{exactRN}) and Lemma \ref{gcd}, the action of $\G$ on $\bD$ is of type $\tqq$, if $X$ is bipartite and of type $\tq$ otherwise.

In order to complete the proof of Theorem \ref{main2}, it only remains to prove that the factor $\Li(\bD)\rtimes \G$ is hyperfinite.  
By \cite{zi}, this follows from the next result.

\begin{proposition}\label{amenable} The action of $\G$ on $\bD$ is amenable.
\end{proposition}

\begin{proof}
The group $G=\aut(\D)$ acts transitively on $\bD$  \cite [Chapter I.8]{ftn}.
Fix an element  $\w \in \bD$, and let
$G_{\w} = \{g \in G : g\w = \w \}$.
Then $\bD \cong G/G_{\w}$, and $\mu$ corresponds to a measure in the unique quasi-invariant measure class of $G/G_{\w}$.
The group  $G_{\w}$ is amenable by \cite[Theorem 8.3]{ftn}.
It follows from \cite[Corollary 4.3.7]{zim} that the action of $\G$ on $G/G_{\w}$ is amenable.
\end{proof}

\bigskip

\section{Appendix: Proof of a Technical Lemma}\label{technicallemma}

This section contains a proof of the technical result, Lemma \ref{bcondition}.
During the course of the proof it will be necessary to concatenate paths.
A difficulty arises because two proper paths cannot necessarily be concatenated to produce a proper path. The product  path may backtrack at the initial edge of the second path. This problem is overcome by introducing a detour around a proper cycle attached at the initial vertex of the second path. The following auxiliary Lemma will be used to do this.

\begin {lemma}\label{aloop} {\sc(Attaching a Loop to an Edge.)}
Let $X$ be a finite connected graph whose vertices all have degree at least three.
Let $e$ be an edge of $X$. Then there is a proper cycle $L$ based at the terminal vertex $t(e)$, not passing through $e$ and having length $\ell(L)\le \delta + \lambda$, where $\delta$ is the diameter of $X$ and $\lambda$ is the maximum length of a circuit in $X$.
\end{lemma}

\begin{proof}
The edge $e$ is contained in a maximal tree $T$ in $X$. Every vertex of $X$ is a vertex of $T$.  Let $P$ be a maximal proper (geodesic) path in $T$ with initial vertex $o(e)$ and initial edge $e$. Let $v$ be the terminal vertex of $P$ and $f$ the terminal edge.
Then $v$ is an endpoint of $T$. The vertex $v$ has degree at least three. It follows that there are two edges in $X$ other than $\overline f$ with initial vertex $v$.
These two edges may both have terminal vertex $v$ (in fact one may be the opposite of the other) or else one or both of them may end at a vertex other than $v$. However in all cases we may use one or both of these edges together with edges in $T$ to construct a circuit $L_0$ based at $v$ and not passing through $e$.
The required proper cycle $L$ can the be constructed from $P\cup L_0$.
\end{proof}

\noindent {\em Proof of {\rm Lemma \ref{bcondition}}.}
Let $\delta$ be the diameter of $X$ and $\lambda$ the maximum length of a circuit in $X$.
We show that property $\fL(K)$ is satisfied with $K=10+ 10\delta +6\lambda$.

Let $P_1, P_2$ be proper paths in $X$ having the same length $n$ and the same initial vertex $v_0$. Let $p_1, p_2$ be the terminal vertices of $P_1, P_2$ respectively. We must construct proper cycles $C_1, C_2$ based at $v_0$ satisfying the conditions of Definition \ref{fL}.

Choose once and for all a path $[p_1,p_2]$ of shortest length between $p_1$ and $p_2$.
There are two separate cases to consider.

\smallskip

\noindent {\sc Case 1}. The length of $[p_1,p_2]$ is even. Denote this length by $2s$ where $s\ge 0$ and let $p_0$ be the midpoint of $[p_1,p_2]$.  If $s=0$ then a simpler argument will apply, and produce a smaller bound for the lengths of $C_1, C_2$, so we assume that $s>0$.

\begin{figure}[htbp]
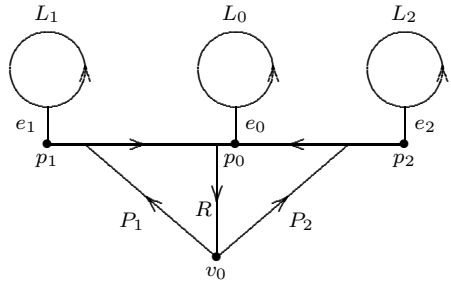
\label{technical}
\hfil
\centerline{
\beginpicture 
\setcoordinatesystem units <0.5cm, 0.5cm>  
\setplotarea  x from -7 to 7,  y from -3 to 3.3
\putrule from -4.5 0 to 5 0 
\putrule from -4.5 0 to -4.5 1
\putrule from 0.5 0  to  0.5 1
\putrule from 5 0 to 5 1
\putrule from 0 0 to 0 -3
\setlinear
\plot -3.5 0  0 -3   3.5 0 /
\put{$_{L_0}$}[b]  at  0.5  3.25 
\put{$_{L_1}$}[b]  at  -4.5  3.25
\put{$_{L_2}$}[b]  at  5  3.25
\put{$_{e_1}$}[r]  at  -4.75  0.5
\put{$_{e_0}$}[l]  at  0.75  0.5
\put{$_{e_2}$}[l]  at  5.25  0.5
\put{$_{P_1}$}[t,r]  at  -1.9  -1.8
\put{$_{P_2}$}[t,l]  at  1.9  -1.8
\put{$_{R}$}[r]  at  -0.1  -1.7
\put{$_\bullet$}  at  -4.5 0
\put{$_{p_1}$}[t]  at  -4.5 -0.25
\put{$_\bullet$}  at  5 0
\put{$_{p_2}$}[t]  at  5 -0.25
\put{$_\bullet$}  at  0 -3
\put{$_{v_0}$}[t]  at  0 -3.25
\put{$_\bullet$}  at  0.5 0
\put{$_{p_0}$}[t]  at  0.5 -0.25
\circulararc 360 degrees from -4.5 1 center at -4.5 2 
\circulararc 360 degrees from 0.5 1 center at 0.5 2
\circulararc 360 degrees from 5 1 center at 5 2
\arrow <6pt> [.3,.67] from  1.75 -1.5 to  1.866 -1.4
\arrow <6pt> [.3,.67] from  -1.75 -1.5 to  -1.866 -1.4
\arrow <6pt> [.3,.67] from  0 -1.4 to  0 -1.5
\arrow <6pt> [.3,.67] from   -2 0 to  -1.9 0
\arrow <6pt> [.3,.67] from   2 0 to  1.9 0
\arrow <6pt> [.3,.67] from  -3.5 2 to  -3.5  2.1
\arrow <6pt> [.3,.67] from  1.5 2 to  1.5  2.1
\arrow <6pt> [.3,.67] from  6 2 to  6  2.1
\endpicture
}
\hfil
\caption{Constructing proper cycles of the same length.}
\end{figure}

Choose a path $R$ of minimal length from $p_0$ to $v_0$. The cycles $C_1, C_2$ will be constructed from portions the paths $P_1, P_2, R, [p_1,p_2]$, with loops attached to avoid backtracking. Refer to Figure \ref{technical}.

Choose an edge $e_1$ with $o(e_1)=p_1$ such that $e_1$ is not the initial edge of $[p_1,p_2]$ and $\ove_1$ is not the final edge of $P_1$.
Choose an edge $e_2$ with $o(e_2)=p_2$ such that $e_2$ is not the initial edge of $[p_2,p_1]$ and $\ove_2$ is not the final edge of $P_2$.
For $i=1,2$, attach a proper cycle $L_i$ at $t(e_i)$, as in Lemma \ref{aloop}.

Assume that the initial edge of $R$ does not meet either of the edges of $[p_1,p_2]$ which contain $p_0$. Let $C_1$ be the proper cycle based at $v_0$ obtained by passing through the following sequence of paths and edges in the order indicated.
\[
P_1\to e_1\to L_1\to \ove_1\to [p_1,p_2]\to e_2\to L_2\to \ove_2 \to [p_2,p_0] \to R
\]
Similarly, let $C_2$ be obtained from
\[
P_2\to e_2\to L_2\to \ove_2\to [p_2,p_1]\to e_1\to L_1\to \ove_1 \to [p_1,p_0] \to R
\] 
The proper cycles $C_1, C_2$ have initial segments $P_1, P_2$ respectively and have the same length $n+k$, where $k=4+\ell(L_1)+\ell(L_2)+3s +\ell(R)\le 4 +2(\delta+\lambda)+
\frac{3}{2}\delta +\delta<4+5\delta+ 2\lambda$.

Now assume that the initial edge of $R$  meets an edge of $[p_1,p_2]$ which contains $p_0$.  This is precisely the situation illustrated in Figure \ref{technical}. The cycles $C_1, C_2$ described above will no longer both be proper, since there will be a backtrack for one of them at the first edge of $R$. In order to avoid this, choose an edge $e_0$ with $o(e_0)=p_0$ such that $e_0$ does not meet either of the edges of $[p_1,p_2]$ containing $p_0$. Attach a proper cycle $L_0$ at $t(e_0)$, as in Lemma \ref{aloop}. Modify the cycles $C_1, C_2$ above so that the final part of each becomes 
\[
\dots, p_0]\to e_0\to L_0\to \ove_0 \to R
\]
The proper cycles $C_1, C_2$ now have the same length $n+k$, where $k< 6+6\delta+ 3\lambda$.

\refstepcounter{picture}

\smallskip 

\noindent {\sc Case 2}. The length of $[p_1,p_2]$ is odd. Denote this length by $2s+1$ where $s\ge 0$ and let $p_0$ be the vertex of $[p_1,p_2]$ with $d(p_0,p_1)=s, d(p_0,p_2)=s+1$. The argument that follows will be slightly different, but simpler, if $s=0$. We therefore again assume that $s>0$. 

Exactly the same argument as in Case 1 shows that there are proper cycles $C_1, C_2$ based at $v_0$ and with initial segments $P_1, P_2$ respectively. The only difference is that
$\ell(C_1)=n+k+1$, $\ell(C_2)=n+k$, where $k< 6+6\delta+ 3\lambda$.

The cycles will be modified to have the same length by adding to the end of each an appropriate proper cycle based at $v_0$.  
The (possibly improper) cycle
\[
P_1\to [p_1,p_2] \to \overline P_2
\]
has odd length. Deleting appropriate parts of this cycle shows that $X$ contains a circuit $C_0$ of odd length $2t+1$. (In other words, the graph $X$ is not bipartite.)

Choose a path $S_1$ of minimal length from $v_0$ to $C_0$. Let $v_1$ be the terminal vertex of $S_1$. The circuit $C_0$ is the union of two proper paths $C_0^+$, $C_0^-$
with lengths $t+1$, $t$ respectively and initial vertex $v_1$.
Let $v_2$ be the terminal vertex of the paths  $C_0^+$, $C_0^-$.
Choose a path $S_2$ of minimal length from $v_2$ to $v_1$.
Add to the end of each of the cycles $C_1, C_2$ a cycle based at $v_0$, as indicated below
\[
C_1\to S_1 \to C_0^-\to S_2
\]
\[
C_2\to S_1 \to C_0^+\to S_2
\]
The resulting cycles have the same length, namely $n+k+1+t+\ell(S_1)+\ell(S_2)=n+k'$, where $k'\le k+\lambda+2\delta< 6+ 8\delta+ 4\lambda$. Either or both of these cycles may have backtracking at $v_0$ or at $v_2$ (but not at $v_1$). If this happens add an edge (and its reverse) and adjoin a loop to both cycles at the relevant vertex as in Lemma \ref{aloop}.
The resulting cycles are proper (i.e. have no backtracking) and have the same length 
$n+k''$, where $k''\le 10+ 10\delta +6\lambda$. 
\qed


\begin{thebibliography}{RRR}

\bibitem [AD]{ad} C. Anantharaman-Delaroche, $C^*$-alg\`ebres de Cuntz-Krieger et groupes Fuchsiens, {\it Operator Theory, Operator Algebras and Related Topics (Timi\c{s}oara 1996)}, 17--35, The Theta Foundation, Bucharest, 1997.

\bibitem[C1]{c} J. Cuntz, A class of $C^*$-algebras and topological Markov chains:
Reducible chains and the Ext-functor for $C^*$-algebras, {\it Invent. Math.} {\bf 63} (1981), 23--50.

\bibitem[C2]{c'} J. Cuntz, K-theory for certain $C^*$-algebras, {\it Ann. of Math.} {\bf 113} (1981), 181--197. 

\bibitem[CK]{ck} J. Cuntz and W. Krieger, A class of $C^*$-algebras and topological Markov chains, {\it Invent. Math.} {\bf 56} (1980), 251--268.

\bibitem[FTN]{ftn} A.~Fig\`a-Talamanca and C.~Nebbia, {\em Harmonic Analysis and Representation Theory for Groups Acting on Homogeneous Trees}, LMS Lecture Note Series, 182, Cambridge University Press, 1991. 

\bibitem [GH]{gh} E. Ghys and P. de la Harpe (editors), {\it Sur les Groupes Hyperboliques
d'apr\`es Mikhael Gromov}, Birkh\"auser, Basel, 1990.

\bibitem[Gr]{gr} M. Gromov, Hyperbolic groups. {\it Essays in group theory}, 75--263, Math. Sci. Res. Inst. Publ., {\bf 8}, Springer, New York, 1987.

\bibitem[HO]{ho} T.~Hamachi and M.~Osikawa, {\em Ergodic Groups of
Automorphisms and Krieger's Theorems}, Seminar on  Mathematical Sciences
No.~3, Keio University, Japan, 1981.

\bibitem [HN]{hn}H. Moriyoshi and T. Natsume, The Godbillon-Vey cyclic cocycle and longitudinal Dirac operators, {\it Pacific J. Math.} {\bf 172} (1996), 483--539.

\bibitem [K]{k} E. Kirchberg, Exact $C^*$-algebras, tensor products, and the classification of purely infinite algebras, {\it Proceedings of the International Congress of Mathematicians (Z\"urich, 1994)}, Vol.~2, 943--954, Birkh\"auser, Basel, 1995.

\bibitem [N]{natsume} T. Natsume,  Euler characteristic and the class of unit in $K$-theory, {\it  Math. Z.} {\bf 194} (1987), 237--243.

\bibitem[Ok]{ok} R. Okayasu, 
Type III factors arising from Cuntz-Krieger algebras, 
{\it Proc. Amer. Math. Soc.} {\bf 131} (2003), 2145--2153.
 
\bibitem[Ped]{ped} G. K. Pedersen, {\it $C^*$-algebras and their Automorphism Groups},
Academic Press, New York, 1979.

\bibitem[Ph]{ph} N. C. Phillips, A classification theorem for nuclear purely infinite simple $C^*$-algebras, {\it Doc. Math.}{\bf 5} (2000), 49--114. 

\bibitem[RR]{rr} J.~Ramagge and G.~Robertson, Factors from trees,  {\it Proc. Amer. Math. Soc. } {\bf 125} (1997), 2051--2055.

\bibitem[Rob]{rob} G. Robertson Torsion in $K$-theory for boundary actions on affine buildings of type $\tilde A_n$, {\it K-theory} {\bf 22} (2001), 251--269.

\bibitem[RS]{rs2} G. Robertson and T. Steger, Affine buildings, tiling systems and higher rank Cuntz-Krieger algebras, {\it J. reine angew. Math.} {\bf 513} (1999), 115--144.

\bibitem[S]{ser} J.-P. Serre, {\it Arbres, amalgames, $SL_2$}, $3^e$ ed., Ast\'erisque {\bf 46}, Soc. Math. France, 1983.

\bibitem[Spa]{spa} R.~J.~Spatzier,  An example of an amenable action from
geometry, {\it Ergod.\ Th.\ \& Dynam.\ Sys.} {\bf 7} (1987), 289--293.

\bibitem[Spi]{sp} J. Spielberg, Free product groups, Cuntz-Krieger algebras, and covariant maps, {\it International J. Math.} {\bf 2} (1991), 457-476.

\bibitem[Su]{su} V.~S.~Sunder, {\em An Invitation to von Neumann Algebras},
Universitext, Springer-Verlag, New York 1987.

\bibitem[Z1]{zi} R.~L.~Zimmer, Hyperfinite factors and amenable ergodic actions, 
{\it Invent. Math.} {\bf 41} (1977), 23--31.

\bibitem[Z2]{zim} R.~L.~Zimmer, {\it Ergodic Theory and Semisimple Groups},
Birkha\"{u}ser, Boston 1985.

\end{thebibliography}
\end{document}